\documentclass[article,11pt,plain]{aiaa-pretty}    % journal, submit, article
\footerlabel{\textit{Oper. Res. Forum} \textbf{4}, 31 (2023) https://doi.org/10.1007/s43069-023-00215-6}

\usepackage[utf8]{inputenc}  % Helps a lot with arXiv submission!!!
\newtheorem{theorem}{Theorem}%  meant for continuous numbers
\newtheorem{proposition}[theorem]{Proposition}%

\newtheorem{proof}{Proof}

\newtheorem{remark}{Remark}%

\newtheorem{definition}{Definition}%

\newtheorem{lemma}{Lemma}
\newcommand{\Real}{\mathbb R}
\newcommand{\myset}[1]{\left\{#1\right\}}
\newcommand{\real}[1]{{\mathbb R}^{#1}}

\newcommand{\bq}{{\boldsymbol q}}

\newcommand{\bu}{{\boldsymbol u}}

\newcommand{\bx}{\boldsymbol x}

\newcommand{\buf}{{\bu(\cdot)}}  % f for function
\newcommand{\bQ}{{\boldsymbol Q}}

\newcommand{\bW}{{\boldsymbol W}}

\newcommand{\U}{\mathbb{U}}

\newcommand{\bzero}{{\bf 0}}

\newcommand{\blam}{{\mbox{\boldmath $\lambda$}}}

\author{ %
I. M. Ross\thanks{Distinguished Professor and Program Director, Control and Optimization} 
\\
\textit{Naval Postgraduate School, Monterey, CA 93943}
}

\title{Derivation of Coordinate Descent Algorithms from Optimal Control Theory}

\abstract{
Recently, it was posited that disparate optimization algorithms may be coalesced in terms of a central source emanating from optimal control theory.
Here we further this proposition by showing how coordinate descent algorithms may be derived from this emerging new principle.
In particular, we show that basic coordinate descent algorithms
can be derived using a maximum principle and a collection of  max functions as ``control'' Lyapunov functions.
The convergence of the resulting coordinate descent algorithms is thus connected to the controlled dissipation of their corresponding Lyapunov functions.
The operational metric for the search vector in all cases
is given by the Hessian of the convex objective function.

}

\begin{document}
\maketitle

\section{Introduction}

A new theory for optimization over the field of real numbers was posited in \cite{rossJCAM-1}.  A proposition in \cite{rossJCAM-1} is that all optimization algorithms can be deduced from some, universal, high-level principle.  The presumptive high-level principle is optimal control theory\cite{vinter}. Admittedly\cite{rossJCAM-1}, the proposed idea appears manifestly ill-conceived because optimal control problems are typically solved using mathematical programming techniques and not the other way around\cite{bazaraa-2006}. Nonetheless, motivated largely by intellectual curiosity and the draw of a unified theory, it was shown in \cite{rossJCAM-1} that the optimal control approach was indeed intrinsically capable of providing a unified framework for optimization. A key detail in this theory is the use of a control Lyapunov function (CLF)\cite{CLF-sontag-sussmann-96}. A CLF for optimization is similar to a merit function in the sense that it measures progress of the iterates toward a solution; however, unlike a merit function, a CLF must be chosen a priori to the construction of an algorithm.  That is, the details of the presumptive algorithm are unknown until additional aspects of controlling the dissipation of the Lyapunov function are worked out.  This process is in sharp contrast to the the more conventional approach of constructing (non-control) Lyapunov functions a posteriori to the design of algorithms for the purposes of proving convergence\cite{polyak-lyap}. In other words, the optimal-control/CLF framework comes with an initial guarantee of convergence of an algorithm, wherein the details of the presumptive algorithm are, at first, unknown but deduced by subsequent steps.  Using this unconventional approach, it was shown in \cite{rossJCAM-1} that dissipating smooth quadratic CLFs produced the standard gradient and Newton's method.  The distinction between these two methods, according to optimal control theory, is the controlled dissipation rate of the CLF.
Nonsmooth Lyapunov functions were not considered in \cite{rossJCAM-1}.  Here we consider a family of elementary nonsmooth Lyapunov functions given by certain max functions and show that they generate coordinate descent (CD) algorithms\cite{wright:CD-survey}.  The proof that CD algorithms may also be subsumed under the unified theory of optimal control is the main contribution of this paper.
The repercussions of the fact that a CD algorithm is part of a larger, new theory implies certain benefits of membership, the details of which are described shortly.

Prior to the advent of the machine learning (ML) revolution\cite{MLopt:SIREV}, CD algorithms were not serious contenders to optimization methods\cite{wright:CD-survey,nesteov-CD-2012}.  The widespread and practical applications of ML to such problems as image recognition, cancer detection, self-driving automobiles and other big-data problems\cite{bengio:deepLearning} changed prevailing attitudes about previously-discarded optimization algorithms\cite{wright:CD-survey,MLopt:SIREV,nesteov-CD-2012}.  CD algorithms are now more widely used than ever before in such diverse application as medical imaging\cite{CD:med-image}, internet-of-things\cite{CD-IoT}, computational biology\cite{CD:bio}, and much more\cite{wright:CD-survey,MLopt:SIREV}. By absorbing it under the emerging unified theory for optimization, we anticipate new formalisms for those ML algorithms that utilize CD methods.  For instance, by linking a CD algorithm to a (nonsmooth) CLF, the meaning of descent can be better stated in terms of generating a decrement of the associated Lyapunov function.  We also show that the operational metric for CD algorithms is hidden and given by the Hessian of the (convex) objective function.  That is, even though the Hessian is not explicitly used in CD algorithms, its properties are crucial for convergence under the rubric of optimal control theory.

\section{Background: Optimal Control Theory for Optimization}
As a matter of completeness, we briefly review certain results from \cite{rossJCAM-1} that are pertinent to a derivation of CD algorithms.  To this end, consider the unconstrained optimization problem given by,
\begin{equation}\label{eq:prob-S}
(S) \left\{\displaystyle\mathop\text{Minimize }_{\bx_f \in \real{N_x}}   E(\bx_f) \right.
\end{equation}
where, $E : \real{N_x} \to \Real, \ N_x \in \mathbb{N}^+$ is a twice continuously differentiable function.  Let, $\bx^0 \in \real{N_x}$ be an initial point of an algorithm designed to solve Problem~$(S)$.  That is, the said algorithm generates a sequence of points,
\begin{equation}\label{eq:algorithm=p2s-map}
\bx^0 \mapsto \myset{\bx^0=\bx_0, \bx_1, \ldots, \bx_k, \bx_{k+1}, \ldots } \quad \bx_k \in \real{N_x}, \ k \in \mathbb{N}
\end{equation}
Let,
\begin{equation}\label{eq:y-def}
y(t) := E(\bx(t))
\end{equation}
where, $\Real \ni t \mapsto \bx(t) \in \real{N_x}$ is a continuously differentiable ``trajectory.''  Suppose that a trajectory passes through the point $\bx^0$ at some point in ``time'' $t = t_0$ such that $\bx(t_0) = \bx^0$. Then, a suitable discretization of this trajectory generates the point-to-set map,
\begin{equation}\label{eq:algorithm-discretized}
\bx^0 \mapsto \myset{\bx^0=\bx(t_0), \bx(t_1), \ldots, \bx(t_k), \bx(t_{k+1}), \ldots } \quad t_{k+1} > t_k \in \Real
\end{equation}
Equation~\eqref{eq:algorithm-discretized} matches \eqref{eq:algorithm=p2s-map} if $\bx(t_k) = \bx_k, \forall\ k \in \mathbb{N}$.
%--------------------
\begin{definition}[Algorithm Primitive]\label{def:algor-primitive}
Given a point $\bx^0 \in \real{N_x}$, the point-to-set map,
\begin{equation}\label{eq:algorithm-primitive-def}
\bx^0 \mapsto \myset{\bx^0 = \bx(t_0), [t_0, \infty) \ni t \mapsto \bx(t)}
\end{equation}
is called an algorithm primitive for solving Problem~$(S)$.
\end{definition}
%---------------------
Defintition~\ref{def:algor-primitive} implies that we consider the search for a suitable algorithm primitive to be a fundamental problem for solving Problem~$(S)$.  To design this fundamental problem, we differentiate \eqref{eq:y-def} with respect to $t$; this yields,
\begin{equation}\label{eq:ydot=}
\dot y(t)= \partial_{\bx} E(\bx(t))\cdot \dot\bx(t)
\end{equation}
Let,
\begin{equation}\label{eq:xdot=u}
\dot\bx(t) = \bu(t)
\end{equation}
where, $t \mapsto \bu(t) \in \real{N_x}$ is a control trajectory (yet to be determined).  Substituting \eqref{eq:xdot=u} in \eqref{eq:ydot=} we generate the following time-free optimal control problem,
\begin{eqnarray}
&(R) \left\{
\begin{array}{lrl}
\textsf{Minimize }  & J[y(\cdot), \bx(\cdot), \buf, t_f]
:=& y(t_f)  \\
\textsf{Subject to} & \dot\bx(t)=& \bu(t)\\
&\dot y(t)=& \partial_{\bx} E(\bx(t))\cdot \bu(t)\\
&\bx(t_0) =& \bx^0 \\
& y(t_0) = & E(\bx^0)
\end{array} \right.& \label{eq:prob-R}
\end{eqnarray}
where, $y(\cdot), \bx(\cdot)$ and $\buf$ are unknown functions, $t_f \ge t_0$ is an unknown final time (parameter) and $J : \big(y(\cdot), \bx(\cdot), \buf, t_f \big) \mapsto \Real$ is the cost functional given by the final-time value of $y(t)$. Because the optimal final-time value of $y(t)$ is the minimum of Problem~$(R)$, it follows that:
%-------------------
\begin{enumerate}
\item Minimizing the cost functional of Problem~$(R)$ is the same as minimizing the objective function of Problem~$(S)$;
\item An optimal system trajectory $[t_0, t_f] \mapsto \big(y(t), \bx(t), \bu(t)\big)$ that solves Problem~$(R)$ contains $[t_0, t_f] \mapsto \bx(t)$.  Thus, an algorithm primitive for solving Problem~$(S)$ is a natural outcome of solving Problem~$(R)$; and,
\item The final-time value of $\bx(t)$ given by $\bx(t_f)$ is a solution to Problem~$(S)$.
\end{enumerate}
%-------------------
The Pontryagin Hamiltonian\cite{ross-book} for Problem~$(R)$ is given by,
\begin{equation}
H(\blam_x, \lambda_y, \bx, y, \bu) := \blam_x \cdot \bu + \lambda_y\, \partial_{\bx}E(\bx)\cdot\bu
\end{equation}
where, the pair $(\bx, \bu)$ is the instantaneous value of $(\bx(t), \bu(t))$ and $(\blam_x, \lambda_y)$ is the adjoint covector pair (costate) that satisfies the (adjoint) differential equations,
\begin{subequations}\label{eq:adj}
\begin{align}
\dot\blam_{x}(t) & := -\partial_{\bx} H = -\lambda_y(t)\,\partial^2_{\bx}E(\bx(t))\, \bu(t) \label{eq:adj-x-unc}\\
\dot\lambda_y(t) & := -\partial_y H = 0 \label{eq:adj-y-unc}
\end{align}
\end{subequations}
The terminal transversality conditions\cite{ross-book} for Problem~$(R)$ are given by
\begin{subequations}\label{eq:TVC}
\begin{align}
\blam_{x}(t_f) & = \bzero \label{eq:TVC-1}\\
\lambda_y(t_f) & = \nu_0 > 0  \label{eq:TVC-2}
\end{align}
\end{subequations}
where \eqref{eq:TVC-2} incorporates the nontriviality condition\cite{vinter,ross-book}.  Substituting \eqref{eq:TVC-2} and \eqref{eq:xdot=u} in \eqref{eq:adj-x-unc} generates the following integral of motion:
\begin{equation}\label{eq:lamx=-gradE}
\blam_x(t) = -\nu_0\, \partial_{\bx} E(\bx(t))
\end{equation}
%-------------
See \cite{rossJCAM-1} for additional details associated with \eqref{eq:lamx=-gradE}.
%===================================
\begin{theorem}[Transversality Mapping\cite{rossJCAM-1}]\label{thm:tmt}
The terminal transversality conditions for Problem $(R)$ generate Fermat's optimality condition for Problem~$(S)$.  Furthermore, all extremals of Problem~$(R)$ are zero-Hamiltonian, infinite order, normal singular arcs.
\end{theorem}
%----------------------------------
\begin{proof}
From \eqref{eq:lamx=-gradE} and \eqref{eq:TVC} it follows that,
\begin{equation}
\partial_{\bx} E(\bx(t_f)) = \bzero
\end{equation}
which is the Fermat condition for Problem~$(S)$ for $\bx_f = \bx(t_f)$.  See \cite{rossJCAM-1} for a proof of the remainder of the statements of this theorem.
\end{proof}
%======================================
\begin{remark}
Problem~$(R)$ was first formulated by Goh\cite{Goh-1997}; however, because the problem is singular; i.e., an application of Pontryagin's Principle does not generate a candidate optimal control, Goh\cite{Goh-1997} and Lee et al\cite{Goh-2021} pursued an alternative theory based on bang-bang controls by adding control constraints to Problem~$(R)$. See also \cite{rossJCAM-1} for a constrained-optimization version of Theorem~\ref{thm:tmt}.
\end{remark}
%------------------------------------

\section{A Deductive Approach to Designing Algorithm Primitives and Algorithms}

From Section~2, it follows that any candidate controller that satisfies the adjoint equations and the transversality conditions, namely, \eqref{eq:adj} and \eqref{eq:TVC}, produces an algorithm primitive (cf.~\eqref{eq:algorithm-primitive-def}).  Consequently, the road to designing an algorithm primitive for Problem~$(S)$ can be stated in terms of the following fundamental problem: given a point $\bx^0$, find a control trajectory, $\bu(t; \bx^0)$ that satisfies the following equation:
\begin{equation}\label{eq:aux-system}
\dot\blam_x(t) = -\nu_0\,\partial^2_{\bx} E(\bx(t))\, \bu(t),    \qquad \blam_x(t_f) = \bzero, \qquad \nu_0 > 0
\end{equation}
The governing equation for the algorithm primitive is then given by (cf.~\eqref{eq:xdot=u}),
\begin{equation}\label{eq:xdot=u-v2}
\dot\bx(t) = \bu(t; \bx^0)
\end{equation}
An Euler discretization of \eqref{eq:xdot=u-v2} generates an algorithm given by,
\begin{equation}\label{eq:Euler-disc}
\bx_{k+1} = \bx_k + h_k \bu_k
\end{equation}
where, $\bx_k = \bx(t_k), \bx(t_0) = \bx^0, \bu_k = \bu(t_k; \bx^0)$ and $h_k > 0$ is the Eulerian step size. To illustrate this idea, consider a candidate control trajectory, $\bu^a(\cdot)$, given by the formula,
\begin{equation}\label{eq:ua:=}
\bu^a(\blam_x(t), \bx(t)) := \big[\partial^2_{\bx} E(\bx(t))\big]^{-1}\, \blam_x(t)
\end{equation}
where we assume the Hessian is invertible.  Substituting \eqref{eq:ua:=} in \eqref{eq:aux-system} we get
\begin{equation}\label{eq:dotLam4Newton}
\dot\blam_x(t) = -\nu_0\,\blam_x(t),   \qquad \nu_0 > 0
\end{equation}
Hence, it follows that $\blam_x(t_0)$ will be driven to zero by \eqref{eq:ua:=}. In other words, \eqref{eq:ua:=} together with \eqref{eq:xdot=u} and \eqref{eq:Euler-disc} generate an algorithm primitive and an algorithm for solving Problem~$(S)$.  To better appreciate this last statement, substitute \eqref{eq:ua:=} in \eqref{eq:Euler-disc} by setting $\bu_k = \bu^a(\blam_x(t_k), \bx(t_k))$; this generates,
\begin{equation}\label{eq:Newton-redux}
\bx_{k+1} = \bx_k - h_k \nu_0\, \big[\partial^2_{\bx} E(\bx_k)\big]^{-1}\, \partial_{\bx} E(\bx_k)
\end{equation}
where, we have eliminated $\blam_x(t)$ according to the integral of motion given by \eqref{eq:lamx=-gradE}.  Obviously, \eqref{eq:Newton-redux} is a ``reinvention'' of Newton's method (damped) with a step size of $\alpha_k := h_k \nu_0$. Note, however, that we did not start with a Newton's method; rather, it was ``derived'' from optimal control theory.

The main point of the illustrative process of the previous paragraph was to show how an algorithm can be derived using \eqref{eq:aux-system} as a foundational equation for Problem~$(S)$.  Except for the first step, i.e., the production of \eqref{eq:ua:=}, the rest of the steps are procedural. In fact, it is apparent that \eqref{eq:ua:=} was ``derived'' based on an anticipation of \eqref{eq:dotLam4Newton}.  In a quest to formalize the production of of control functions like \eqref{eq:ua:=} that are anticipatory of convergence, we use the notion of control Lyapunov functions (CLFs) that were briefly introduced in Section~1.  In simple terms, a CLF for \eqref{eq:aux-system} is a positive-definite function $V : \blam_x \to \Real_+$ which dissipates with an application of some $\bu$ whenever the value of $V$ is nonzero. See \cite{CLF-sontag-sussmann-96,clarkeLyap} for technical details. The motivation for a CLF in control theory is largely driven by its applications to stabilizing physical systems.  As a result, the ideas of stability, controllability and stabilizability are interwoven with the use of a CLF.  For optimization purposes, a CLF may be used in a less restrictive sense; hence, we forgo the usual technicalities and note that the dissipation criterion (for a smooth CLF) requires the existence of a control $\bu$ such that $d_tV(\blam_x(t))  < 0 $ whenever $\blam_x(t) \ne \bzero$.  This simple idea translates to,
\begin{equation}\label{eq:Vdot=}
 \partial_{{\footnotesize{\blam}}_x}V(\blam_x) \cdot \dot\blam_x =  -\langle \partial_{{\footnotesize{\blam}}_x}V(\blam_x),\ \nu_0\,\partial^2_{\bx} E(\bx)\, \bu \rangle < 0
\end{equation}
whenever $\blam_x \ne \bzero$.  Thus, for a maximal dissipation rate, we choose $\bu$ to maximize the inner product given given by $\langle \cdot, \cdot \rangle $ in \eqref{eq:Vdot=}. This idea is formalized in terms of the following \emph{maximum principle}:
\begin{eqnarray}
&(P) \left\{
\begin{array} {lll}
\displaystyle\mathop\textsf{Maximize }_{\bu \in \U({\footnotesize{\blam}}_x, \bx, y, t)}  &\langle \partial_{{\footnotesize{\blam}}_x} V(\blam_x), \partial^2_{\bx} E(\bx)\, \bu \rangle
\end{array} \right.& \label{eq:Max-P}
\end{eqnarray}
where, $\U(\blam_x, \bx, y, t)$ is any compact set in $\real{N_x}$ that may depend upon $(\blam_x, \bx, y, t)$.  If a CLF is nonsmooth (but continuous), then from the results of \cite{clarkeLyap}, \eqref{eq:Vdot=} must be modified to hold for all $ \partial_{{\footnotesize{\blam}}_x} V(\blam_x) \in \partial V(\blam_x)$,
where $\partial_{{\footnotesize{\blam}}_x} V(\blam_x)$ is a subgradient of  $V(\blam_x)$ and $\partial V(\blam_x)$ is a subdifferential of $V$.  Thus, \eqref{eq:Max-P} generalizes to,
\begin{eqnarray}
&(P') \left\{
\begin{array} {lll}
\displaystyle\mathop\textsf{Maximize }_{\bu \in \U({\footnotesize{\blam}}_x, \bx, y, t)}  &\langle \partial_{{\footnotesize{\blam}}_x} V(\blam_x), \partial^2_{\bx} E(\bx)\, \bu \rangle
\end{array} \right.& \quad \forall\  \partial_{{\footnotesize{\blam}}_x} V(\blam_x) \in \partial V(\blam_x) \label{eq:Max-P-nonsmooth}
\end{eqnarray}
The maximum principle, given by \eqref{eq:Max-P-nonsmooth}, forms part of the new foundations for designing algorithm primitives for Problem~$(S)$. In terms of this foundation, an algorithm primitive is determined a posteriori to the selection of $V$ and $\U$. The maximum principle also formalizes and generalizes certain other ideas that have been pursued in earlier works; see, for example \cite{bhaya-book}, where well-established results from control theory have been used quite successfully in solving certain optimization problems.

\section{Selection of $(V, \U)$ Pairs That Generate Basic Coordinate Descent Algorithms}
Because the objective function in \eqref{eq:Max-P} is linear in the control variable, we choose $\U$ in a manner that metricizes the underlying space\cite{rossJCAM-1}; i.e., we set,
\begin{equation}\label{eq:U=quad}
\U(\blam_x, \bx, y, t) := \myset{\bu:\ \langle \bu,\ \bW(\blam_x, \bx, y, t)\, \bu \rangle \le \Delta(\blam_x, \bx, y, t) }
\end{equation}
where $\bW: (\blam_x, \bx, y, t) \mapsto \mathbb{S}^{N_x}_{++}$ is a positive definite matrix function and $\Delta : (\blam_x, \bx, y, t) \mapsto \Real_{++}$ is a strictly positive function whose value resembles a trust region in optimization\cite{NW:NumOptBook}.
%==============================================

%==================================
\begin{lemma}\label{lemma:1}
Let, $V: \blam_x \mapsto \Real$ be a continuous CLF such that $\partial^2_{\bx}E(\bx)\,\partial_{{\footnotesize{\blam}}_x} V(\blam_x) \ne \bzero$ for all $\partial_{{\footnotesize{\blam}}_x} V(\blam_x) \in \partial V(\blam_x) $ and $\blam_x \ne \bzero$.  Let $\U$ be given by \eqref{eq:U=quad}.  Then, the maximum principle generates a control $\bu$ that satisfies the dissipation condition given by \eqref{eq:Vdot=} for all $\partial_{{\footnotesize{\blam}}_x} V(\blam_x) \in \partial V(\blam_x) $ and $\blam_x \ne \bzero$.
%dissipates $V(\blam_x)$ for \eqref{eq:aux-system}.
\end{lemma}
%-------------------------------------
\begin{proof}
Select some element $ \partial_{{\footnotesize{\blam}}_x} V(\blam_x) \in \partial V(\blam_x)$.  An application of the maximum principle given by \eqref{eq:Max-P-nonsmooth} with $\U$ stipulated by \eqref{eq:U=quad} generates the control law,
\begin{equation}\label{eq:Wu=4quadU}
\bW(\blam_x, \bx, y, t)\, \bu = \sigma[@t]\, \partial^2_{\bx} E(\bx)\,\partial_{{\footnotesize{\blam}}_x} V(\blam_x)
\end{equation}
where, $\sigma[@t] \equiv \sigma(\blam_x(t), \bx(t), y, t) > 0$ is a positive scalar function. Substituting \eqref{eq:Wu=4quadU} in \eqref{eq:Vdot=} we get,
\begin{multline}
 \langle \partial_{{\footnotesize{\blam}}_x}V(\blam_x),\ \partial^2_{\bx} E(\bx)\, \bu \rangle \\
=  \sigma[@t] \ \Big[ \partial^2_{\bx} E(\bx)\, \partial_{{\footnotesize{\blam}}_x}V(\blam_x)\Big]^T \bW^{-1}(\blam_x, \bx, y, t)\Big[ \partial^2_{\bx} E(\bx)\, \partial_{{\footnotesize{\blam}}_x}V(\blam_x)\Big] > 0
\end{multline}
for all $\partial_{{\footnotesize{\blam}}_x} V(\blam_x) \in \partial V(\blam_x) $ and $\blam_x \ne \bzero$.
%whenever $\blam_x \ne \bzero$.
\end{proof}
%=====================================

Now, consider a candidate CLF given by,
\begin{equation}\label{eq:CLF=maxQuad}
V^C(\blam_x) := \max_{1 \le i \le N_x} \myset{V^C_i(\blam_x)}, \quad V^C_i(\blam_x) := (\lambda_x^i)^2/2
\end{equation}
where, $\lambda_x^i, i = 1, \ldots, N_x$ are the components of $\blam_x$.  Then $\partial V^C(\blam_x)$ is formally given by\cite{HUL-book-96,BL-book-2006}
\begin{equation}\label{eq:subdiff-formula}
\partial V^C(\blam_x) := \mathbf{co} \cup \myset{\partial V^C_i(\blam_x):\ V^C_i(\blam_x) = V^C(\blam_x)}
\end{equation}
The subdifferential $\partial V^C_i(\blam_x)$ in \eqref{eq:subdiff-formula} is the same as the gradient, $\partial_{{\footnotesize{\blam}}_x}V^C_i(\blam_x)$; namely,
\begin{equation}
\partial V^C_i(\blam_x) = \partial_{{\footnotesize{\blam}}_x}V^C_i(\blam_x) = (0, \ldots, \lambda_x^i, \ldots, 0)
\end{equation}
That is, the components of $\partial_{{\footnotesize{\blam}}_x}V^C_i(\blam_x)$ are all zero except for the $i^{th}$ element whose value is given by $\lambda_x^i$. Furthermore, for the case $\left(\lambda^i_x\right)^2 \ne \left(\lambda^j_x\right)^2, i \ne j = 1, \ldots, N_x$, \eqref{eq:subdiff-formula} is a singleton and can be written as,
\begin{equation}
\partial V^C(\blam_x)  = \myset{ \bq \circ \blam_x } = \partial_{{\footnotesize{\blam}}_x}V^C(\blam_x)
\end{equation}
where, $\circ$ denotes the Hadamard product and $\bq \in \real{N_x}$ is a vector whose components $q^i, \ i = 1, \ldots, N_x$ are defined by
\begin{equation}\label{eq:q=def}
q^i = \left\{
        \begin{array}{lll}
          1, & \hbox{if } &\left(\lambda^i_x\right)^2 \in \displaystyle{\max_j}\myset{ \left(\lambda^j_x\right)^2, \ j = 1, \ldots, N_x } \\
          0, & \hbox{if } &otherwise
        \end{array}
      \right.
\end{equation}
In \eqref{eq:q=def}, the dependence of $q^i$ on $\blam_x$ is suppressed for notational convenience.

In all other cases (i.e., when $\left(\lambda^i_x\right)^2 = \left(\lambda^j_x\right)^2, i \ne j$)
$\partial V^C(\blam_x)$ is not a singleton and given by,
\begin{equation}\label{eq:V-subdiff-explicit}
\partial V^C(\blam_x) := \mathbf{co} \myset{(0, \ldots, \lambda_x^i, \ldots, 0):\ \left(\lambda^j_x\right)^2/2 = V^C(\blam_x)}
\end{equation}
Furthermore, \eqref{eq:q=def} generates a vector $\bq$ with multiple components taking the value unity. As a result, $\bq \circ \blam_x$ is along a subgradient of $V^C(\blam_x)$ that is equally weighted in terms of the gradients of $V_i^C(\blam_x)$.
In this situation, we may use \eqref{eq:q=def} to construct an ``unbiased'' subgradient of $V^C(\blam_x)$  as,
\begin{equation}\label{eq:gradVC=gamma*}
\partial_{{\footnotesize{\blam}}_x}V^C(\blam_x) = \gamma\, \bq \circ \blam_x, \quad \gamma > 0
\end{equation}
where, $\gamma$ is chosen to satisfy the condition, $\gamma\, \bq \circ \blam_x \in \partial V^C(\blam_x)$.
%================================
\begin{proposition}\label{prop:deriv-CDstd}
Suppose $E$ is a strictly convex, twice-differentiable function.  Let, $\bW(\bx) := \partial^2_{\bx} E(\bx)$ and $\U \equiv \U(\bx)$ be given by \eqref{eq:U=quad} Then, the maximum principle with $V(\blam_x) = V^C(\blam_x)$  generates a candidate coordinate descent algorithm given by,
\begin{equation}\label{eq:CD-std}
\bx_{k+1} = \bx_k - \alpha_k\, \bq \circ \partial_{\bx} E(\bx_k)
\end{equation}
where, $\alpha_k > 0$ is the algorithm step size.
\end{proposition}
%--------------------------------
\begin{proof}
From Lemma~\ref{lemma:1}, it follows that for $\bW(\bx) = \partial^2_{\bx} E(\bx)$ and  $V(\blam_x) = V^C(\blam_x)$ the maximum principle generates the control law,
\begin{equation}\label{eq:u4stdCD}
\bu = \sigma[@t]\, \partial_{{\footnotesize{\blam}}_x} V^C(\blam_x) = \sigma[@t]\, \gamma\, \bq \circ \blam_x
\end{equation}
where, the second equality in \eqref{eq:u4stdCD} uses the unbiased subgradient of $V^C(\blam_x)$ given by \eqref{eq:gradVC=gamma*}.  Substituting \eqref{eq:u4stdCD} in \eqref{eq:xdot=u} and using the singular-arc integral of motion given by \eqref{eq:lamx=-gradE} we get,
\begin{equation}\label{eq:xdot=stdCD}
\dot\bx = - \sigma[@t]\,\gamma\, \nu_0\, \bq \circ \partial_{\bx} E(\bx)
\end{equation}
A forward Euler discretization of \eqref{eq:xdot=stdCD} with a step size of $h_k > 0$ generates \eqref{eq:CD-std} with $\alpha_k := h_k\, \sigma_k\, \gamma_k\, \nu_0 $, where, $\sigma_k = \sigma[@t_k]$.
\end{proof}
%==================================
\begin{remark}
It is apparent from the proof of Proposition~\ref{prop:deriv-CDstd} that the Eulerian step size, $h_k$, is not the same as the optimization step size, $\alpha_k$.
\end{remark}
%---------------------------------

In ML applications, $N_x$ is very large.  In these instances, it is desirable to produce search directions for blocks of variables that may be parsed for parallel computation.  Motivated by this need, we partition $\blam_x$ into $m$-blocks of $n_i, i = 1, \ldots, m$ variables that allows us to write,
\begin{equation}
\blam_x = (\lambda^1_x, \lambda^2_x, \ldots, \lambda^{N_x}_x) = (\blam^{n_1}_x, \ldots, \blam^{n_{m}}_x)
\end{equation}
where, $\sum_{i=1}^m n_i = N_x$.  Define,
\begin{subequations}
\begin{align}
V^B_i\big(\blam_x\big) &= \frac{1}{2}\langle\blam^{n_i}_x,\, \bQ^i \blam^{n_i}_x \rangle, \quad \bQ^i \in \mathbb{S}^{n_i}_{++}\\
%V^\Sigma(\blam_x) &= \sum_{i=1}^m V^i\big(\blam^{n_i}_x\big) \\
V^B(\blam_x) &= \max_{1\le i \le m} \myset{ V^B_i\big(\blam_x\big) }
\end{align}
\end{subequations}
Then, using arguments similar to those that generated \eqref{eq:gradVC=gamma*} we get,
\begin{equation}
\partial_{{\footnotesize{\blam}}_x} V^B(\blam_x) = \gamma\, \bq \circ \bQ\,\blam_x, \quad \gamma > 0
\end{equation}
where, $\bQ = \text{diag}(\bQ^i, i = 1, \ldots, m)$ and $\bq$ is appropriately modified in its definition in \eqref{eq:q=def} with $\left(\lambda^i_x\right)^2$ replaced by $V^B_i(\blam_x)$ and $1$ replaced by a vector of $n_i$-ones.
%
%============================
\begin{proposition}
Let the assumptions of Proposition~\ref{prop:deriv-CDstd} hold.  Then, the maximum principle with $V(\blam_x) = V^B(\blam_x)$ generates a block coordinate descent algorithm given by,
\begin{equation}\label{eq:blockCD}
\bx_{k+1} = \bx_k - \alpha_k\, \bq \circ \bQ\, \partial_{\bx} E(\bx_k)
\end{equation}
where, $\alpha_k > 0$ is the algorithm step size.
\end{proposition}
%----------------------------
\begin{proof}
The proof is almost identical to that of Proposition~\ref{prop:deriv-CDstd}; hence, it is omitted for brevity.
\end{proof}
%======================================

Finally, we consider yet another candidate CLF given by,
\begin{equation}\label{eq:CLF=infNorm}
V^\infty(\blam_x) := \big\|\blam_x\big\|_{\infty}
\end{equation}
Rewriting \eqref{eq:CLF=infNorm} as
\begin{equation}
V^\infty(\blam_x) := \max_{1 \le i \le N_x} \myset{V^\infty_i(\blam_x)}, \quad V^\infty_i(\blam_x) := \abs{\lambda_x^i}
\end{equation}
it follows from \eqref{eq:subdiff-formula} through \eqref{eq:gradVC=gamma*} that we may select the unbiased subgradient of $V^\infty(\blam_x)$ according to,
\begin{equation}
\partial_{{\footnotesize{\blam}}_x}  V^\infty(\blam_x) = \gamma\,\bq \circ \text{sign}(\blam_x), \quad \gamma > 0
\end{equation}
with the definition of $\bq$ modified in \eqref{eq:q=def} with $\left(\lambda^i_x\right)^2$ replaced by $\abs{\lambda^i_x}$.
%
%============================
\begin{proposition}
Let the assumptions of Proposition~\ref{prop:deriv-CDstd} hold.  Then, the maximum principle with $V(\blam_x) = V^\infty(\blam_x)$ generates a normalized coordinate descent algorithm given by,
\begin{equation}\label{eq:CD-normalized}
\bx_{k+1} = \bx_k - \alpha_k\, \bq \circ \textrm{sign}\big(\partial_{\bx} E(\bx_k)\big)
\end{equation}
where, $\alpha_k > 0$ is the algorithm step size.
\end{proposition}
%---------------------------
\begin{proof}
The steps in this proof are nearly identical to that of Proposition~\ref{prop:deriv-CDstd}; hence, they are omitted for brevity.
\end{proof}
%======================================
\begin{remark}
Equation~\eqref{eq:CD-normalized} is generated in \cite{boyd-book} by seeking a steepest descent method in the $\ell_1$ norm. The optimal control route to producing \eqref{eq:CD-normalized} provides the insight that the $\ell_\infty$ norm of the gradient of the (convex) objective function is a CLF for \eqref{eq:aux-system}.
\end{remark}
%-------------------------------
\begin{remark}
In Propositions 2, 3 and 4, the control space is metricized by the Hessian.  The control variable in all these instances is a continuous-time proxy for the search vector in optimization. In this context, \eqref{eq:U=quad} may be interpreted as a generalized- and control-version of the variable metric concept originally proposed by Davidon\cite{davidon} in 1959.
\end{remark}

%=====================================================================

%\backmatter

%\bmhead{Acknowledgments}
%
%Funding for this research was provided by the the Air Force Office of Scientific Research.

\section*{Declarations}

\begin{itemize}
\item My manuscript has no associated data
\item On behalf of all authors, the corresponding author states that there is no conflict of interest.
\item Funding for this research was provided by the the Air Force Office of Scientific Research.
\end{itemize}



\begin{thebibliography}{10}

\bibitem{rossJCAM-1}
I. M. Ross, An optimal control theory for nonlinear optimization, J. Comp. and Appl. Math., 354, 2019, 39--51.
doi 10.1016/j.cam.2018.12.044.

\bibitem{vinter}
R. B. Vinter,  Optimal Control, Birkh\"{a}user, Boston, MA,
2000.

\bibitem{bazaraa-2006}
M. S. Bazaraa, H. D. Sherali, C. M. Shetty, Nonlinear Programming: Theory and Algorithms, 3rd Edition, Wiley-Interscience, Hoboken, NJ, 2006, 367--369.


\bibitem{CLF-sontag-sussmann-96}
E. D. Sontag,  H. J. Sussmann, General classes of control-Lyapunov functions, ISNM International Series of Numerical Mathematics, 121, Birkhauser Basel, 1996, 87--97.
%In: Jeltsch, R., Mansour, M. (eds) Stability Theory. https://doi.org/10.1007/978-3-0348-9208-7_10


\bibitem{polyak-lyap}
B. Polyak, P. Shcherbakov, Lyapunov functions: an optimization theory perspective, IFAC PapersOnLine, 50-1, 2017, 7456--7461.

\bibitem{wright:CD-survey}
S. J. Wright, Coordinate descent algorithms, Mathematical Programming, 151/1, 2015, 3--34

\bibitem{MLopt:SIREV}
L. Bottou, F. E. Curtis, J. Nocedal, Optimization methods for large-scale machine learning, SIAM Review, 2018 60/2, 223--311.


\bibitem{nesteov-CD-2012}
Yu. Nesterov, Efficiency of coordinate descent methods on huge-scale optimization problems, SIAM J. Opt., 22/2, 2012, 341--362


\bibitem{bengio:deepLearning}
Y. Bengio, Y. LeCun, G. Hinton, Deep Learning, Nature, 521, 2015, 436--444.

\bibitem{CD:med-image}
J.-B. Thibault, K. D. Sauer, C. A. Bouman, J. Hsieh, A three-dimensional statistical approach to improved image quality for multislice helical CT, Med. Phys., 34/11, 2007, 4526--4544.

\bibitem{CD-IoT}
H. Wang et al., A comprehensive survey on training acceleration for large machine learning models in IoT, IEEE Internet of Things Journal, 9/2, 2022, 939--963.

\bibitem{CD:bio}
G. Ifrim, W. Carsten, Bounded coordinate-descent for biological sequence classification in high dimensional predictor space, Proc. 17th ACM SIGKDD Intl. Conf. on knowledge discovery and data mining, KDD 2011, 708--716.


\bibitem{ross-book}
I. M. Ross, A Primer on Pontryagin's Principle in Optimal Control, second ed., Collegiate Publishers, San Francisco, CA, 2015.


\bibitem{Goh-1997}
B. S. Goh, Algorithms for unconstrained optimization via control theory, J. Optim. Theory Appl., 92/3, 1997, pp.~581--604.

\bibitem{Goh-2021}
M. S. Lee, H. G. Harno, B. S. Goh, K. H. Lim, On the bang-bang control approach via a component-wise line search strategy for unconstrained optimization, Numerical Algebra, Control and Optimization, 11/1, 2021, pp.~45--61.

\bibitem{clarkeLyap}
F. Clarke, Lyapunov functions and feedback in nonlinear control. In: M.S. de Queiroz, M. Malisoff, P. Wolenski (eds) Optimal control, stabilization and nonsmooth analysis. Lecture Notes in Control and Information Science, vol 301. Springer, Berlin, Heidelberg (2004), 267--282.


\bibitem{bhaya-book}
A. Bhaya, E. Kaszkurewicz, Control Perspectives on Numerical Algorithms and Matrix Problems, SIAM, Philadelphia, PA, 2006.


\bibitem{NW:NumOptBook}
J. Nocedal, S. Wright, Numerical Optimization, Second Ed, Springer, New York, NY, 2006.


\bibitem{HUL-book-96}
J-B. Hiriart-Urruty and C. Lemar\'{e}chal, Convex Analysis and Minimization Algorithms I, Springer-Verlag, Berlin Heidelberg GmbH, 1996.

\bibitem{BL-book-2006}
J. Borwein, A. Lewis, Convex Analysis and Nonlinear Optimization: Theory and Examples, Second Ed, Springer, New York, N.Y., 2006.


\bibitem{boyd-book}
S. Boyd, L. Vandenberghe, Convex Optimization, Cambridge University Press, Cambridge, UK, 2009.

\bibitem{davidon}
W. C. Davidon, Variable metric method for minimization, SIAM J. optimization, 1/1 (1991), pp.~1--17 (originally published as Argonne National Laboratory Research and Development Report 5990, May 1959; revised November 1959).



\end{thebibliography}


\end{document}